\documentclass[11pt]{article}

\usepackage[margin=1in]{geometry}
\usepackage[T1]{fontenc}
\usepackage{lmodern}
\usepackage{microtype}
\usepackage{amsmath,amssymb,amsthm,mathtools}
\usepackage{graphicx}
\usepackage{subcaption}
\usepackage{comment}
\usepackage{tikz}
\usepackage[colorlinks=true,linkcolor=blue,citecolor=blue,urlcolor=cyan]{hyperref}

\hypersetup{
    pdftitle={A Classification of Trees Without Unique Subgraphs},
    pdfauthor={Andrew Searns},
    pdfkeywords={graph reconstruction, anchor-free graph, unique subgraph, tree},
}

\tikzset{
    vertex/.style = {
        circle,
        draw=black,
        thick,
        minimum size=3mm,
        inner sep=0pt
    },
    every edge/.style = {
        line width=0.8pt
    }
}

\usetikzlibrary{shapes.misc}
\usetikzlibrary{shapes.arrows}
\usetikzlibrary{shapes.geometric}

\newcommand{\N}{\mathbb{N}}
\newcommand{\Z}{\mathbb{Z}}

\DeclareMathOperator*{\argmin}{arg\,min}
\DeclareMathOperator{\rad}{rad}
\DeclareMathOperator{\diam}{diam}
\DeclareMathOperator{\aut}{Aut}

\theoremstyle{plain}
\newtheorem{theorem}{Theorem}
\newtheorem{lemma}{Lemma}[section]
\newtheorem{corollary}[lemma]{Corollary}
\newtheorem{proposition}[lemma]{Proposition}
\theoremstyle{definition}
\newtheorem{definition}{Definition}[section]

\title{A Classification of Trees Without Unique Subgraphs}
\author{Andrew Searns\\[2pt]
    \small Johns Hopkins University\\
    \small\texttt{asearns1@jh.edu}}
\date{\today}

\begin{document}

\maketitle

\begin{abstract}
    We study trees without unique subgraphs, in both the induced and non-induced settings. Graphs without unique induced subgraphs, called anchor-free graphs, form a natural candidate class for counterexamples to the Graph Reconstruction Conjecture, but their structure remains poorly understood. We develop a framework for studying anchor-free graphs based on a covering lemma that characterizes anchor-freeness via induced subgraph coverings of vertices. Using this framework, we prove that a tree is anchor-free if and only if it is a path with one less vertex than a prime. We also introduce the related notion of unique-subgraph-freeness, which acts as an edge-based analogue of anchor-freeness. We show that the line graph of every unique-subgraph-free graph is anchor-free and obtain a complete classification of unique-subgraph-free trees as either prime paths, stars, or extended stars. Together, these results provide a classification of trees without unique subgraphs and suggest a general approach for extending the theory beyond trees.
\end{abstract}

\noindent\textbf{Keywords:} Graph Reconstruction, Anchor, Anchor-Free Graph, Unique Subgraph, Tree, Induced Subgraph, Automorphism Orbit

\smallskip

\section{Introduction}
The Graph Reconstruction Conjecture has been described as the most beautiful and important open conjecture in the field of graph theory \cite{harary1974survey}. Originating from the work of Kelly and Ulam \cite{kelly1942isometric,ulam1960collection}, the conjecture asks how much information is needed to uniquely reconstruct a graph from its parts. Today, even one of the most straightforward and simplest cases of the problem remains unsolved: can we reconstruct a graph with at least three vertices from its deck of vertex deleted subgraphs?

Research on the vertex reconstruction conjecture primarily comes in two directions: probabilistic techniques aiming prove the conjecture generally but possibly missing specific counterexamples, and constructive techniques building up increasingly large classes of graphs which are reconstructible. These dual approaches offer distinct strengths and weaknesses which provide insight into the difficulty in proving the conjecture as a whole. In this paper we hope to bridge some of the gap between the two approaches by studying an explicit family of graphs not covered by probabilistic approaches.

In a seminal probabilistic paper, Bollob\'as proved that almost every graph is reconstructible \cite{bollobas1990almost}. The core theorem of the proof shows that almost every graph has a unique, asymmetric subgraph on $n-2$ vertices. Intuitively given such a subgraph, the uniqueness guarantees that it can be identified in both cards containing it. By asymmetry of the subgraph, the missing edges incident to it in one card can be reconstructed by referencing the opposite card. The only question that remains is whether the two additional vertices are connected by an edge. This is easily answered using combinatorial tools like Kelly's Lemma which constructively counts the number of each proper subgraph from the deck of a graph \cite{kelly1957congruence}.

Bollob\'as's approach suggests that the hardest class of graphs to reconstruct are those with large symmetry groups. The remaining graphs are those which do not have any unique, asymmetric $n-2$ vertex subgraphs. A specific sub-family of these graphs was recently studied by Farhadian who coined the term anchor to describe a unique induced subgraph \cite{farhadian2017simple}. In her paper, Farhadian introduced tools for working with anchors and anchor-free graphs. She used these tools to give an intuitive proof of the reconstruction conjecture for small graphs with at most nine vertices, in contrast with the enumerative verification of the reconstruction conjecture for at most 11 vertices by McKay \cite{mckay1997small, mckay2022reconstruction}.

While much attention has been dedicated to the study of unique subgraphs, to our knowledge, Farhadian's work is the only such work dedicated to the study of anchor-free graphs. Most of the previous work in unique subgraphs relates to extremal type problems, often trying to bound the maximum numbers of unique subgraphs that can occur in a graph. 

For an upper bound, \cite{bradavc2025unique} proved that the maximum possible ratio of distinct subgraphs out of the approximately $2^{\binom{n}{2}}/n!$ subgraphs contained in any $n$ vertex graph approaches $0$ as $n \rightarrow \infty$. On the lower bound, there exist graphs with at least $2^{n^2/2 - cn^{3/2}}$ unique subgraphs \cite{entringer1972number}.

When considering graphs with few distinct subgraphs, it is known that if $G$ has at most $\epsilon n^2$ distinct subgraphs, then $G$ has a large independent set or clique \cite{alon1989graphs,erdHos1989number}. A stronger converse is also known to be true: if $G$ has no clique or independent set on $c_1 \log(n)$ vertices, then $G$ has at least $2^{c_2 n}$ distinct induced subgraphs \cite{shelah1998erdos}.

We also observe that the problem of finding graphs with minimal unique subgraphs can also be viewed through a different dual perspective where the goal is to maximize the number of times a particular subgraph occurs as an induced subgraph. These packing problems have been well studied (e.g. see \cite{yap1988packing} for a survey of results in this field).

Despite all of these results on the structure of unique subgraphs, little is known about graphs without any unique (induced) subgraphs. In this paper, we extend the tools of Farhadian to provide more insights on anchor-free graphs. Our main contribution is a characterization of anchor-free trees as paths with one less vertex than a prime. We also introduce a new notion of unique-subgraph-freeness which is subtly different from anchor-freeness. Many of the techniques for anchor-freeness carry over to unique-subgraph-freeness, and we highlight the important distinctions between the two notions.

\section{Preliminaries}

\subsection{Notation}
All graphs in this paper are finite, simple, and undirected. A graph $H = (V', E')$ is a subgraph of $G = (V, E)$ if $V' \subseteq V$ and $E' \subseteq E \cap (V' \times V')$. If $V' \subsetneq V$, then $H$ is a proper subgraph. We say that $H$ is an induced subgraph of $G$ if $E' = E \cap (V' \times V')$, and in this case we write $H = G[V']$. When it is clear from context, we sometimes abuse the notation and write $G[H]$ to denote $G[V(H)]$.

An isomorphism between graphs $G$ and $H$ is a bijection $\phi : V(G) \to V(H)$ that preserves adjacency. An isomorphism from $G$ to itself is called an automorphism, and the set of all automorphisms of $G$ forms a group under composition, denoted $\aut(G)$. The action of $\aut(G)$ on $V(G)$ partitions the vertices into orbits. Two vertices $u$ and $v$ are in the same orbit if there exists $\phi \in \aut(G)$ with $\phi(u) = v$; in this case we write $u \sim v$. Similarly, the action of $\aut(G)$ on $E(G)$ partitions the edges into edge orbits. 

For graphs $G$ and $H$, we write $s(G, H)$ for the number of subgraphs of $G$ isomorphic to $H$, and $i(G, H)$ for the number of induced subgraphs of $G$ isomorphic to $H$.

\subsection{Unique Subgraphs and Anchors}

An induced subgraph $H$ is an \textit{anchor} for $G$ if $H$ occurs uniquely among induced subgraphs of $G$. That is $H$ is an anchor for $G$ if and only if $i(G, H) = 1$. A graph $G$ is \textit{anchor-free} if its only anchors are the trivial induced subgraphs: the empty graph and $G$ itself.

A subgraph $H \subseteq G$ is a \textit{unique subgraph} if $s(G, H) = 1$. A graph $G$ is \textit{unique-subgraph-free} if the only unique subgraphs of $G$ are the trivial subgraphs: the empty graphs on $0$ and $|V(G)|$ vertices and $G$ itself. \\

We will later show that anchor-freeness and unique-subgraph-freeness are distinguishable properties: there exist graphs with each property which lack the other. Despite this, there are important similarities between the two properties which makes it worthwhile to consider both in tandem. For example, the following proposition shows that any unique subgraph necessarily induces an anchor.

\begin{proposition}\label{prop:unique_proper_subgraphs_anchors}
    If $H$ is a unique subgraph of $G$, then the induced subgraph $G[H]$ is an anchor.
\end{proposition}
\begin{proof}
    Let $A = G[H]$ and suppose that $A$ is not an anchor. Then there is another induced subgraph $B \subseteq G$ such that $A \cong B$ but $B \neq A$. In particular, let $\phi: A \rightarrow B$ be an isomorphism. Now since $A$ and $B$ are different induced subgraphs, they have different vertex sets. However, since $H$ is a spanning subgraph of $A$, the subgraph $\phi(H)$ is a spanning subgraph of $B$. Now $H$ and $\phi(H)$ have different vertex sets, contradicting that $H$ was a unique subgraph.
\end{proof}

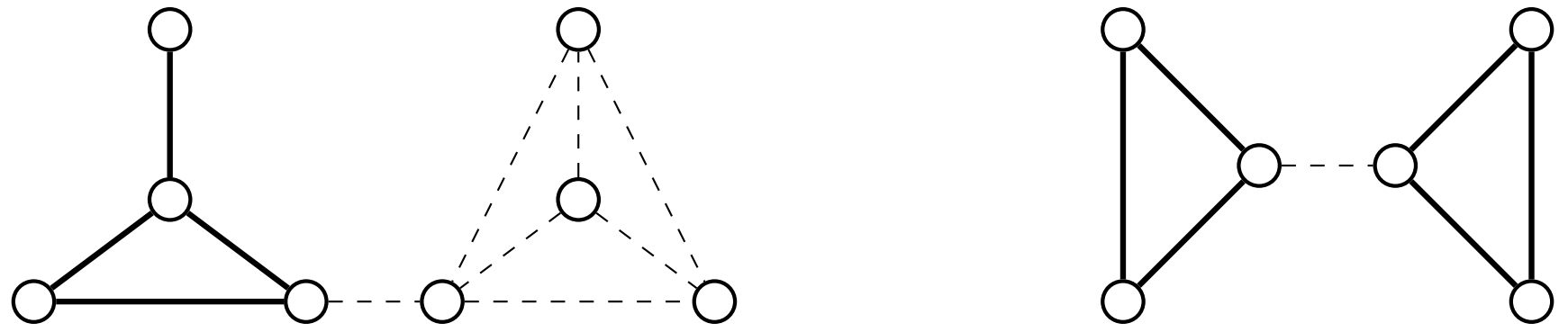
\begin{figure}
    \centering
    \begin{tikzpicture}  
        \node[vertex] (p1)      at (0, -1/4) {};
        \node[vertex] (p2)      at (0, 1) {};
        \node[vertex] (p3)      at (-1, -1) {};
        \node[vertex] (p4)      at (1, -1) {};
        \node[vertex] (k1)      at (3, -1/4) {};
        \node[vertex] (k2)      at (3, 1) {};
        \node[vertex] (k3)      at (2, -1) {};
        \node[vertex] (k4)      at (4, -1) {};
    
        \draw[very thick] (p1) -- (p2);
        \draw[very thick] (p1) -- (p3);
        \draw[very thick] (p1) -- (p4);
        \draw[very thick] (p3) -- (p4);
        \draw[dashed] (p4) -- (k3);
        \draw[dashed] (k1) -- (k2);
        \draw[dashed] (k1) -- (k3);
        \draw[dashed] (k1) -- (k4);
        \draw[dashed] (k2) -- (k3);
        \draw[dashed] (k2) -- (k4);
        \draw[dashed] (k3) -- (k4);
        
        \node[vertex] (a) at (8,0) {};
        \node[vertex] (b) at (7,1) {};
        \node[vertex] (c) at (7,-1) {};
        
        \node[vertex] (d) at (9,0) {};
        \node[vertex] (e) at (10,1) {};
        \node[vertex] (f) at (10,-1) {};
        
        \draw[very thick] (a) -- (b);
        \draw[very thick] (b) -- (c);
        \draw[very thick] (c) -- (a);
        
        \draw[very thick] (d) -- (e);
        \draw[very thick] (e) -- (f);
        \draw[very thick] (f) -- (d);
        
        \draw[dashed] (a) -- (d);
        
    \end{tikzpicture}
    \caption{The difference between anchors and unique subgraphs. The highlighted subgraph on the left is an anchor but not a unique subgraph, and the highlighted subgraph on the right is a unique subgraph but not an anchor.}
    \label{fig:anchor_and_unique_subgraph}
\end{figure}

Proposition~\ref{prop:unique_proper_subgraphs_anchors} shows that unique subgraphs induce anchors. It is tempting to conclude that anchor-freeness implies unique-subgraph-freeness; however this fails when we consider unique subgraphs which span all vertices in a graph (For example, see the second graph of figure~\ref{fig:anchor_and_unique_subgraph}). In this case, a spanning unique subgraph induces the whole graph which is excluded from the definition of anchor-freeness.

\def\L{1.5}
\def\h{0.866*\L}

\section{Properties of Anchor-Free and Unique-Subgraph-Free Graphs}

In \cite{farhadian2017simple}, Farhadian proved many important results for studying anchors and anchor-free graphs. We provide some of these results here as they offer a strong conceptual framework for working with these properties. Where possible, we extend these results to unique subgraphs and unique-subgraph-free graphs.

We start with an important result connecting the study of anchors to the automorphism action of a graph.
\begin{lemma}\label{lem:anchors-respect-orbits}\cite{farhadian2017simple}
    Let $A$ be an anchor for a graph $G$. Let $\mathcal{O}$ be any vertex orbit under the action of $\aut(G)$ on $G$. Then either $A \cap \mathcal{O} = \emptyset$ or $A \cap \mathcal{O} = \mathcal{O}$.
\end{lemma}
\begin{proof}
    Suppose for contradiction that there exists two vertices $u, v \in \mathcal{O}$ with $u \in A$ and $v \notin A$. Since $u \sim v$, there exists an automorphism $\phi \in \aut(G)$ such that $\phi(u) = v$. Now since $\phi$ is an automorphism, $\phi(A) \cong A$. By the choice of $u \in A$, $v = \phi(u) \in \phi(A)$. At the same time, the choice of $v \notin A$ implies that $\phi(A)$ is a different induced subgraph of $G$ isomorphic to $A$. Thus $i(G, A) \geq 2$ contradicting that $A$ is an anchor for $G$.
\end{proof}

A nearly identical argument can be made for unique-subgraph-free graphs. The proof is the same except we replace vertex orbits with edge orbits.

\begin{lemma}\label{lem:usf-respects-orbits}
    Let $A$ be a unique subgraph for a graph $G$. Let $\mathcal{O}$ be any edge orbit under the action of $\aut(G)$ on $G$. Then either $A \cap \mathcal{O} = \emptyset$ or $A \cap \mathcal{O} = \mathcal{O}$.
\end{lemma}
\begin{proof}
    Suppose for contradiction that there exists two edges $e_1, e_2 \in \mathcal{O}$ with $e_1 \in A$ and $e_2 \notin A$. Since $e_1 \sim e_2$, there exists an automorphism $\phi \in \aut(G)$ such that $\phi(e_1) = e_2$. Now since $\phi$ is an automorphism, $\phi(A) \cong A$. By the choice of $e_1 \in A$, $e_2 = \phi(e_1) \in \phi(A)$. At the same time, the choice of $e_2 \notin A$ implies that $\phi(A)$ is a different subgraph of $G$ isomorphic to $A$. Thus $s(G, A) \geq 2$ contradicting that $A$ is a unique subgraph for $G$.
\end{proof}

One important consequence of lemmas~\ref{lem:anchors-respect-orbits}~and~\ref{lem:usf-respects-orbits} is that there are infinitely many anchor-free and unique-subgraph-free graphs.

\begin{corollary}\label{cor:transitive_properties}
    Vertex transitive graphs are anchor-free. Edge-transitive graphs with no isolated vertices are unique-subgraph-free.
\end{corollary}
\begin{proof}
    Vertex transitive graphs have a single vertex orbit. Any anchor must therefore be trivial implying that the graph is anchor-free. A similar argument shows that edge-transitive graphs are unique-subgraph-free.
\end{proof}

A natural next question is whether all anchor-free graphs must be vertex transitive. The smallest non-vertex-transitive, anchor-free graph is $P_4$. By lemma~\ref{lem:anchors-respect-orbits}, any anchor of $P_4$ must respect the orbit partition and thus either consist of a pair of adjacent ($v_2, v_3$) or non-adjacent ($v_1, v_4$) vertices. It is easy to see that neither induced subgraph is an anchor, and thus $P_4$ is anchor-free.

Using an exhaustive search, Farhadian provides a list of all non-vertex-transitive graphs that are anchor-free on at most nine vertices. We note that Farhadian's list erroneously includes $P_8$. We verified the rest of this list using an exhaustive search simplified by lemma~\ref{lem:anchors-respect-orbits} to only count relevant subgraphs. Over the next few lemmas, we build a simple construction to yield infinitely many connected, non-vertex-transitive, anchor-free graphs.

\begin{lemma}\label{lem:af_disjoint_unions}
    Let $E_k$ be the empty graph on $k$ vertices. Then for any connected graph $G$, the cartesian product $G \square E_k$ is anchor-free iff $G$ is anchor-free.
\end{lemma}
\begin{proof}
    Observe that any anchor for $G \square E_k$ must respect the orbit partition of $G \square E_k$. Since $E_k$ is vertex transitive, every orbit $O$ of $G$ is in bijection with the orbits $O \square E_k$ of $G \square E_k$.
    
    Suppose that $G \square E_k$ has an anchor $S$. Since $S$ respects the orbit partition of $G \square E_k$, $S = H \square E_k$ for some $H \subsetneq G$. Observe that if $H$ is not an anchor of $G$ then there is an $H' \subsetneq G$ such that $H' \cong H$ and $H' \square E_k \cong H \square E_k$ contradicting that $S$ is an anchor for $G \square E_k$. Thus $H$ is an anchor for $G$. Therefore if $G \square E_k$ is not anchor-free then $G$ is not anchor-free.

    Now suppose that $G \square E_k$ is anchor-free. Suppose for contradiction that $G$ is not anchor-free. Let $H \subseteq G$ be an anchor for $G$ and consider $H \square E_k$. Since $G \square E_k$ is anchor-free, there is another subgraph $S \subseteq G \square E_k$ with $S \cong H \square E_k$.
    
    Observe first that $H \square E_k$ has $k$ isomorphic copies of $H$ and so $S$ must also have $k$ isomorphic copies of $H$. Since there is only one way to select $H$ from $G$, we conclude that every component of $G \square E_k$ must select exactly one copy of $H$. We conclude that $H \square E_k$ is unique in $G \square E_k$ and thus an anchor. This contradicts that $G \square E_k$ is anchor-free and thus it must be that $G$ is anchor-free.    
\end{proof}

The previous lemma shows that taking arbitrarily many disjoint copies of an anchor-free graph yields an anchor-free graph. In contrast, the next result shows that this structure is the only way that a disconnected graph can be anchor-free.

\begin{lemma}
    Let $G, H$ be connected graphs. Then the disjoint union $G \sqcup H$ is anchor-free iff $G \cong H$ and $G$ is anchor-free.
\end{lemma}
\begin{proof}
    If $G \cong H$ and $G$ is anchor-free, then $G \sqcup H = G \square E_2$ is anchor-free by the previous lemma. Now suppose that $G \ncong H$. If $G$ is not an anchor of $G \sqcup H$, then $G \subseteq H$ as a proper induced subgraph. Thus $H$ is an anchor of $G \sqcup H$.
\end{proof}

The previous lemmas provide a way to build infinitely many non-vertex-transitive, anchor-free graphs: $P_4 \square E_k$ for $k \in \N$.

When we turn our attention back to unique-subgraph-free graphs, the smallest non-edge-transitive, unique-subgraph-free graph is $P_5$. The proof is analogous to the proof that $P_4$ is anchor-free since the action of $\aut(P_5)$ on edges is isomorphic to the action of $\aut(P_4)$ on vertices. This idea can be extended to connect other anchor-free and unique-subgraph-free graphs through line graphs.

\begin{proposition}\label{prop:line_graph_relationship}
    If $G$ is unique-subgraph-free then $L(G)$ is anchor-free.
\end{proposition}
\begin{proof}
    Suppose $G$ is unique-subgraph-free, and let $H$ be a nontrivial induced subgraph of $L(G)$. Define $L^{-1}(H)$ to be the subgraph of $G$ containing each edge $e \in E(G)$ corresponding to a vertex $L(e) \in V(H)$. (We may either take $V(L^{-1}(H))$ to be exactly the vertices incident to the edges of $E(L^{-1}(H))$ or all $V(G)$ arbitrarily.) Since $L^{-1}(H)$ is a subgraph of $G$, there is another subgraph $H'$ of $G$ such that $H' \cong L^{-1}(H)$. By Whitney's Theorem, \nocite{whitney1932congruent} $L(H') \cong L(L^{-1}(H)) = H$ and since $E(H') \neq E(L^{-1}(H))$, $V(L(H')) \neq V(H)$.
    
    It only remains to check that $L(H')$ is an induced subgraph of $L(G)$. Any edge $(uv, vw) \in L(G)$ incident to $L(H')$ corresponds to an edge $(v, w) \in E(G)$ adjacent to an edge $(u, v) \in E(H')$. The line graph edge $(uv, vw)$ is in $L(H')$ exactly when $(v, w) \in E(H')$ and thus $vw \in V(L(H'))$. Thus $L(H')$ is an induced subgraph of $L(G)$.
    
    Since the choice of $H$ was arbitrary, it follows that every induced subgraph of $L(G)$ occurs multiple times, and $L(G)$ is anchor-free.
\end{proof}

While the converse is almost always true, we note that it relies upon Whitney's Theorem which states that if $L(G) \cong L(H)$ then $G \cong H$ with sole exception of $G = K_{3}$ and $H = K_{1,3}$. We can use this exception to build a counterexample to the converse of proposition~\ref{prop:line_graph_relationship}. Let $G = K_{3} \cup K_{1,3}$ which has a unique subgraph of $K_{3}$. Then $L(G) = K_3 \cup K_3$ which is vertex-transitive and thus anchor-free. Thus it is not the case that if $L(G)$ is anchor-free then $G$ is unique-subgraph-free.

Proposition~\ref{prop:line_graph_relationship} highlights the fact that anchor-freeness is a property related to the vertices of a graph while unique-subgraph-freeness is a similar property related to the edges of a graph. With this intuition, one might suspect to find differences between these two properties when looking at complements.

\begin{lemma}\cite{farhadian2017simple}
    If $A$ is an anchor for $G$, then $\overline{A}$ is an anchor for $\overline{G}$.
\end{lemma}
\begin{proof}
    Suppose $G$ has an anchor $A$ but that $\overline{A}$ is not an anchor for $\overline{G}$. Then there is another induced subgraph $B \subseteq \overline{G}$ such that $B \cong \overline{A}$. However, $\overline{B}$ is a proper induced subgraph of $G$ isomorphic to $A$. Since $A$ is an anchor, we then have $\overline{B} = A$ and then $B = \overline{A}$. Since this holds for all $B \cong \overline{A}$, we conclude that $\overline{A}$ is an anchor for $\overline{G}$.
\end{proof}

Since anchors are preserved under taking complements, we see that anchor-freeness is preserved under complements as well.

\begin{corollary}\label{lem:af_complement}
    $G$ is anchor-free iff $\overline{G}$ is.
\end{corollary}

It's worth noting that while anchor-freeness is preserved under taking complements, unique-subgraph-freeness is not. For example, consider the complete bipartite graph $K_{n, m}$ for $n > m$. Since $K_{n, m}$ is edge transitive, it is unique-subgraph-free. $\overline{K_{n, m}} = K_n \cup K_m$ has a unique subgraph in its largest clique $K_n$.

Combining lemma~\ref{lem:af_disjoint_unions} and corollary~\ref{lem:af_complement}, we can show that the join of an anchor-free graph with itself is anchor-free.

\begin{corollary}\label{lem:anchor-free-joins}
    If $G$ is anchor-free the join $G + G$ is anchor-free.
\end{corollary}
\begin{proof}
    $G + G = \overline{\overline{G} \sqcup \overline{G}}$. Two graph $G$ and $\overline{G}$ have complementary subgraphs and thus $G$ is anchor-free iff $\overline{G}$ is anchor-free. The result follows from the previous lemma.
\end{proof}

Using corollary~\ref{lem:anchor-free-joins}, we can now build an infinite sequence of connected, anchor-free, non-vertex-transitive graphs with $G_0 = P_4$, $G_{i+1} \coloneqq G_i + G_i$.

\subsection{Two Useful Covering Lemmas}

While both properties guarantee that subgraphs occur multiple times, neither property provides any information about how isomorphic subgraphs occur in relation to one another. Because of this, it is often more helpful to consider the union of all copies of any subgraph.

\begin{definition}
    Given graphs $G$ and $H$, define the cover of $H$ in $G$ to be the minimal subgraph of $G$ which contains all subgraphs of $G$ isomorphic to $H$. We say $G$ is edge $H$-covered if every edge of $G$ is contained in a subgraph of $G$ isomorphic to $H$.

    Similarly define the induced cover of $H$ in $G$ to be the minimal induced subgraph of $G$ which contains all induced subgraphs of $G$ isomorphic to $H$. We say $G$ is vertex induced $H$-covered if every vertex of $G$ is contained in an induced subgraph of $G$ isomorphic to $H$. 
\end{definition}

Farhadian proved that every anchor-free graph is vertex induced $H$-covered for any induced subgraph $H$ of $G$. We observe that this condition is also sufficient to identify anchor-free graphs.

\begin{lemma}\label{lem:af_covering}
    A graph $G$ is anchor-free iff for every nonempty induced subgraph $H \subseteq G$, $G$ is vertex induced $H$-covered.
\end{lemma}
\begin{proof}
    Suppose first that $G$ has a nontrivial anchor $H$. Since $H$ is a proper subgraph that occurs once as an induced subgraph, $G$ cannot be vertex induced $H$-covered.

    Now suppose for contradiction that $G$ is anchor-free but is not vertex induced $H$-covered for some proper induced subgraph $H$. Let $A$ be the induced cover of $H$ in $G$. Since $G$ is not vertex induced $H$-covered, $A$ is a proper subgraph of $G$. By anchor-freeness of $G$, there is another induced subgraph $B \cong A$. Let $\alpha: A \rightarrow B$ be an isomorphism. Since $V(B) \neq V(A)$ there exists at least one vertex $v \in V(B) - V(A)$. Because $\alpha^{-1}(v) \in A$ which is the induced cover of $H$ in $G$, $\alpha^{-1}(v)$ is contained in an induced subgraph $H'$ of $G$ isomorphic to $H$. Since all of $H'$ is contained in $A$, $\alpha(H') \subseteq B$ is an induced subgraph of $G$ isomorphic to $H$. However, $v \in \alpha(H') - V(A)$ contradicts the maximality of $A$. Thus if $G$ is anchor-free, then for every proper induced subgraph $H$, $G$ is vertex induced $H$-covered.
\end{proof}

We note that in the previous lemma, the choice of using vertex induced coverings was important. For anchor-freeness, it is not possible to extend this lemma to an edge covering. To see this, consider $H = K_3$ in Figure~\ref{fig:anchor_and_unique_subgraph}. In this example, the bridge is not covered by $H$ despite the fact that $G$ is anchor-free.

In contrast, when considering unique-subgraph-free graphs, the natural analog of lemma~\ref{lem:af_covering} uses edge $H$-coverings. The proof is nearly identical. We relax $H$ to be any subgraph, the induced cover to be a cover, and we replace the uncovered vertex $v$ by an uncovered edge $e$.

\begin{lemma}\label{lem:usf_covering}
    A graph $G$ is unique subgraph free iff for every nontrivial subgraph $H$, $G$ is edge $H$-covered.
\end{lemma}

To illustrate the power of lemmas~\ref{lem:af_covering}~and~\ref{lem:usf_covering}, consider a maximum clique. This clique is necessarily an induced subgraph, so in any anchor-free or unique subgraph free graph, every vertex of $G$ must be contained in a clique of maximum size. For anchor-free graphs, lemma~\ref{lem:af_covering} also implies the same result for maximum independent sets. This immediately implies a notable extension of Ramsey's theorem in the context of anchor-free graphs.

\begin{corollary}\label{cor:ramsey}
    If $G$ is anchor-free and $|V(G)| \geq R(s, t)$ then either every vertex of $G$ is contained in a clique of size at least $s$ or every vertex is contained in an independent set of size $t$.
\end{corollary}
\begin{proof}
    Since $|V(G)| > R(s, t)$, Ramsey's theorem requires that $G$ contains either $K_s$ or $I_t$ as an induced subgraph. Lemma~\ref{lem:af_covering} implies that if one of these two induced subgraphs exist, then every vertex is contained in a copy of it.
\end{proof}

In contrast, lemma~\ref{lem:usf_covering} implies a similar result for edge coverings by spanning trees. That is, for each isomorphism class of spanning tree, every edge of a unique-subgraph-free graph must be contained in a spanning tree isomorphic to that that class.

These are perhaps the most promising directions for future work in studying anchor-free and unique-subgraph-free graphs. Corollary~\ref{cor:ramsey} implies strong constraints on any search for large anchor-free graphs.

\section{Anchor-Free Trees}

For this paper, we restrict our attention to anchor-free trees (and forests by lemma~\ref{lem:af_disjoint_unions}). We will eventually prove the following characterization theorem.

\begin{theorem}\label{thm:af_trees}
    The only anchor-free trees are the paths $P_{p-1}$ for $p$ prime.
\end{theorem}
\addtocounter{theorem}{-1}

We divide the proof into two parts. First, we classify the anchor-free paths. We then show that every non-path tree necessarily contains an anchor.

\begin{lemma}\label{lem:anchor-free_paths}
    A path $P_k$ is anchor-free iff $k = p - 1$ for $p$ prime.
\end{lemma}
\begin{proof}
    Since paths are bipartite, any path of odd order has a unique maximum independent set. This independent set is a proper induced subgraph (except in the degenerate case where $k = 1$ where $P_1 = P_{2 - 1}$ is anchor-free) and is therefore a nontrivial anchor.

    Now consider the even path $P_{2k}$. Since anchors respect orbit partitions, any anchor must be invariant under the reflection automorphism of the path. Such a symmetric anchor must arise by deleting the same number of vertices from each side of the center. Suppose without loss of generality that there exists an anchor of $2k - 2d$ vertices. Removing $2d$ vertices from $P_{2k}$ partitions it into $2d + 1$ subpath segments (with some segments possibly containing $0$ vertices). If two segments have different lengths, then we obtain a different isomorphic subgraph by swapping the order of these two segments (See Figure~\ref{fig:path_segment_swap}). Thus any anchor must contain $2d+1$ segments all of the same integer length. That is $\frac{2k - 2d}{2d + 1} = \frac{2k + 1 - 2d - 1}{2d + 1} = \frac{2k + 1}{2d + 1} - 1 \in \Z$.
    
    If for some nontrivial $d$ ($d \neq 0$ and $d \neq k$), $\frac{2k + 1}{2d + 1}$ is an integer, then there is a unique way of removing $2d$ vertices from $P_{2k}$ to yield $2d+1$ segments of equal length and hence form an anchor. Thus the path $P_{2k}$ is anchor-free if and only if there is no $d \neq 0,k$ such that $\frac{2k+1}{2d+1} \in \Z$. Since $2$ does not divide $2k+1$, this condition is equivalent to $2k+1$ being prime. (See Figure~\ref{fig:path_anchors} for an example on $P_{14}$.)
\end{proof}

\begin{figure}
    
        \centering
        \begin{tikzpicture}
        \node[vertex] (a1) at (0,0) {};
        \node[vertex] (a2) at (1,0) {};
        \node[vertex] (a3) at (2,0) {};
        \node[vertex] (a4) at (3,0) {};
        \node[vertex] (a5) at (4,0) {};
        \node[vertex] (a6) at (5,0) {};
        \node[cross out,draw=black,thick,text width=7pt,text height=7pt] at (a2) {};
        \node[cross out,draw=black,thick,text width=7pt,text height=7pt] at (a5) {};
    
        \node[vertex] (b1) at (10,0) {};
        \node[vertex] (b2) at (11,0) {};
        \node[vertex] (b3) at (12,0) {};
        \node[vertex] (b4) at (13,0) {};
        \node[vertex] (b5) at (14,0) {};
        \node[vertex] (b6) at (15,0) {};
        \node[cross out,draw=black,thick,text width=7pt,text height=7pt] at (b3) {};
        \node[cross out,draw=black,thick,text width=7pt,text height=7pt] at (b5) {};
    
        \node[vertex] (c1) at (0,-2) {};
        \node[vertex] (c3) at (2,-2) {};
        \node[vertex] (c4) at (3,-2) {};
        \node[vertex] (c6) at (5,-2) {};
    
        \node[vertex] (d1) at (10,-2) {};
        \node[vertex] (d2) at (11,-2) {};
        \node[vertex] (d4) at (13,-2) {};
        \node[vertex] (d6) at (15,-2) {};
    
        \node at (0,-1) {1};
        \node at (2.5,-1) {2};
        \node at (5,-1) {1};
        \node at (10.5,-1) {2};
        \node at (13,-1) {1};
        \node at (15,-1) {1};
    
        \node[double arrow,draw,minimum height=2cm] at (7.5, 0) {};
        \node at (7.5,-2) {$\cong$};
        
        \draw[dashed] (a1) -- (a2);
        \draw[dashed] (a2) -- (a3);
        \draw (a3) -- (a4);
        \draw[dashed] (a4) -- (a5);
        \draw[dashed] (a5) -- (a6);
    
        \draw (b1) -- (b2);
        \draw[dashed] (b2) -- (b3);
        \draw[dashed] (b3) -- (b4);
        \draw[dashed] (b4) -- (b5);
        \draw[dashed] (b5) -- (b6);
        
        \draw (c3) -- (c4);
    
        \draw (d1) -- (d2);
        
        \end{tikzpicture}
        
        \caption{When the subpath segments have different lengths, reordering them creates isomorphic subgraphs by removing different vertices.}
        
        \label{fig:path_segment_swap}
        \vspace{0.5cm}
    
        \begin{tikzpicture}
        \node[vertex] (a1) at (0,0) {};
        \node[vertex] (a2) at (1,0) {};
        \node[vertex] (a3) at (2,0) {};
        \node[vertex] (a4) at (3,0) {};
        \node[vertex] (a5) at (4,0) {};
        \node[vertex] (a6) at (5,0) {};
        \node[vertex] (a7) at (6,0) {};
        \node[vertex] (a8) at (7,0) {};
        \node[vertex] (a9) at (8,0) {};
        \node[vertex] (a10) at (9,0) {};
        \node[vertex] (a11) at (10,0) {};
        \node[vertex] (a12) at (11,0) {};
        \node[vertex] (a13) at (12,0) {};
        \node[vertex] (a14) at (13,0) {};

        \draw (a1) -- (a2);
        \draw (a2) -- (a3);
        \draw (a3) -- (a4);
        \draw (a4) -- (a5);
        \draw (a5) -- (a6);
        \draw (a6) -- (a7);
        \draw (a7) -- (a8);
        \draw (a8) -- (a9);
        \draw (a9) -- (a10);
        \draw (a10) -- (a11);
        \draw (a11) -- (a12);
        \draw (a12) -- (a13);
        \draw (a13) -- (a14);

        \node[cross out,draw=black,thick,text width=7pt,text height=7pt] at (4,-1) {};
        \node[cross out,draw=black,thick,text width=7pt,text height=7pt] at (9,-1) {};
        \node[vertex] (b1) at (0,-1) {};
        \node[vertex] (b2) at (1,-1) {};
        \node[vertex] (b3) at (2,-1) {};
        \node[vertex] (b4) at (3,-1) {};
        \node[vertex] (b6) at (5,-1) {};
        \node[vertex] (b7) at (6,-1) {};
        \node[vertex] (b8) at (7,-1) {};
        \node[vertex] (b9) at (8,-1) {};
        \node[vertex] (b11) at (10,-1) {};
        \node[vertex] (b12) at (11,-1) {};
        \node[vertex] (b13) at (12,-1) {};
        \node[vertex] (b14) at (13,-1) {};

        \draw (b1) -- (b2);
        \draw (b2) -- (b3);
        \draw (b3) -- (b4);
        \draw (b6) -- (b7);
        \draw (b7) -- (b8);
        \draw (b8) -- (b9);
        \draw (b11) -- (b12);
        \draw (b12) -- (b13);
        \draw (b13) -- (b14);

        \node[cross out,draw=black,thick,text width=7pt,text height=7pt] at (2,-2) {};
        \node[cross out,draw=black,thick,text width=7pt,text height=7pt] at (5,-2) {};
        \node[cross out,draw=black,thick,text width=7pt,text height=7pt] at (8,-2) {};
        \node[cross out,draw=black,thick,text width=7pt,text height=7pt] at (11,-2) {};
        \node[vertex] (c1) at (0,-2) {};
        \node[vertex] (c2) at (1,-2) {};
        \node[vertex] (c4) at (3,-2) {};
        \node[vertex] (c5) at (4,-2) {};
        \node[vertex] (c7) at (6,-2) {};
        \node[vertex] (c8) at (7,-2) {};
        \node[vertex] (c10) at (9,-2) {};
        \node[vertex] (c11) at (10,-2) {};
        \node[vertex] (c13) at (12,-2) {};
        \node[vertex] (c14) at (13,-2) {};

        \draw (c1) -- (c2);
        \draw (c4) -- (c5);
        \draw (c7) -- (c8);
        \draw (c10) -- (c11);
        \draw (c13) -- (c14);
        
        \end{tikzpicture}
        \caption{The anchors of $P_{2k}$ are uniquely determined by the factors of $2k+1$. Here $P_{14}$ has exactly two anchors constructed from $14 + 1 = 3 \times 5 = 5 \times 3$.}
        \label{fig:path_anchors}
\end{figure}

For the remainder of this section, assume that $T$ is an anchor-free, non-path tree. We derive progressively stronger restrictions on the branch structure of $T$, ultimately obtaining a contradiction. We begin by introducing some standard notation which will aid us in our proofs.

The \textit{eccentricity} of a vertex $\epsilon_G(v) = \max_{w \in V(G)} d(v, w)$ is the maximum shortest distance to any other vertex in $G$. When the graph is clear from context, we write $\epsilon(v)$.

The \textit{radius} and \textit{diameter}, denoted $\rad(G), \diam(G)$ respectively, are the minimum and maximum eccentricities of $G$: $\rad(G) = \min_{v \in V(G)} \epsilon(v)$ and $\diam(G) = \max_{v \in V(G)} \epsilon(v)$.

The \textit{center} $Z(G)$ of $G$ is the set of vertices with minimum eccentricity (the radius): $Z(G) = \argmin_{v \in V(G)} \epsilon(v)$. 

Given a tree $T$ and a vertex $v \in V(T)$, we call the components of $T - v$ \textit{branches} of $T$ rooted at the neighbors of $v$. For any vertex $u \in V(T) - v$, we write $B_{v}(u)$ to denote the branch which contains the vertex $u$. We note that $u$ does not need to be in $N(v)$ with this notation and that there may be multiple choices of $u$ which select the same branch. The \textit{height} of a branch $B_v(u)$ is the maximum distance from its root (the neighbor of $v$ in that branch) to any vertex of the branch.

An \textit{end-cutvertex} of a tree is a vertex which is adjacent to a leaf. \\

The following lemmas illustrate the strength of lemma~\ref{lem:af_covering} in the context of trees. By applying it to carefully chosen subgraphs, we obtain strong global restrictions on the branch structure of anchor-free trees.

\begin{lemma}\label{lem:diametral_paths}
    In an anchor-free tree, every vertex lies on a diameter of the graph.
\end{lemma}
\begin{proof}
    This follows directly from lemma~\ref{lem:af_covering} with the subgraph $H = P_{D + 1}$ where $D$ is the diameter of the graph.
\end{proof}

Intuitively, since every vertex lies on a diameter, the leaves of the tree must all be at the same height from the center of the graph.

\begin{lemma}\label{lem:af_leaf_height}
    In an anchor-free tree $T$, every leaf is at the same distance from the center.
\end{lemma}
\begin{proof}
    By lemma~\ref{lem:diametral_paths}, every leaf lies on a diametral path. A leaf lying on a diametral path must be an endpoint of it, and all diametral endpoints are equidistant from the center.
\end{proof}

Any non-path tree must necessarily contain a vertex with degree at least $3$. Anchor-freeness turns the presence of even a single such vertex into restrictions on the degree near leaves of the tree.

\begin{lemma}\label{lem:af_star_covers}
    In an anchor-free tree $T$ with at least one vertex of degree at least $3$, every end-cutvertex has degree at least $3$.
\end{lemma}
\begin{proof}
    Since there exists at least one vertex with degree at least $3$, $K_{1, 3}$ is an induced subgraph of $T$. By lemma~\ref{lem:af_covering}, $T$ is induced $K_{1, 3}$-covered. Each leaf of $T$ has degree $1$ and thus must be adjacent to a vertex of degree at least $3$.
\end{proof}

Lemmas~\ref{lem:af_leaf_height}~and~\ref{lem:af_star_covers} impose strong restrictions on the branch structures of anchor-free trees. In fact, the next lemma implies that if we are given a branch from an anchor-free tree, then we can identify which vertex was the root of the branch.

\begin{lemma}\label{lem:af_root_identification}
    Let $v$ be a vertex of an anchor-free tree $T$ with an adjacent branch $B$ that does not contain a center vertex of $T$. If $B$ is not a star, then from the structure of $B$ alone, we can identify the unique vertex adjacent to $v$ in $T$.
\end{lemma}
\begin{proof}
    Since $B$ is a branch of $T$, lemma~\ref{lem:af_star_covers} implies that in $T$ every end-cutvertex has degree at least $3$. By lemma~\ref{lem:af_leaf_height}, all leaves are at the same height so the end-cutvertices of $T$ cannot be adjacent to a single leaf and two internal neighbors of $T$. Consequently, all leaves of $T$ come in pairs. Thus if $B$ contains an isolated leaf vertex, then this vertex must be $r$.

    Suppose no such distinguishing leaf exists. Then by lemma~\ref{lem:af_leaf_height} every maximal sub-branch of $B$ has the same height. This forces $B$ to have a unique central vertex which must be $r$, since attaching $B$ to $T$ at any other vertex would create branches of unequal heights at $v$. Hence $r$ is uniquely determined.
\end{proof}

It is important to note that in the previous lemma, if $B$ is a star, then since all leaves are in the same orbit of $B$, we can choose any leaf to be the root of $B$. However, this choice is not uniquely determined and will require a small amount of additional casework when we prove the full theorem.

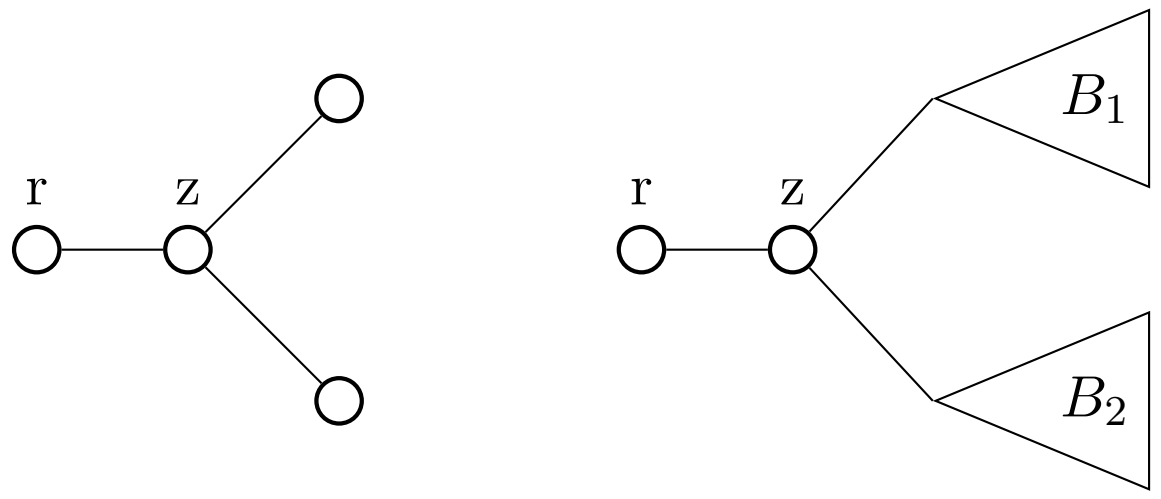
\begin{figure}
    \centering
    \begin{tikzpicture}
            \node[vertex, label=above:r] (a1) at (0,0) {};
            \node[vertex, label=above:z] (a2) at (1,0) {};
            \node[vertex] (a3) at (2,1) {};
            \node[vertex] (a4) at (2,-1) {};

            \draw (a2) -- (a1);
            \draw (a2) -- (a3);
            \draw (a2) -- (a4);

            \node[vertex, label=above:r] (b0) at (4,0) {};
            \node[vertex, label=above:z] (b1) at (5,0) {};
            \node[isosceles triangle,minimum width=1cm,minimum height=1cm,shape border rotate=180,draw] (b2) at (7,1) {$B_1$};
            \node[isosceles triangle,minimum width=1cm,minimum height=1cm,shape border rotate=180,draw] (b3) at (7,-1) {$B_2$};

            \draw (b0) -- (b1);
            \draw (b1) -- (b2.apex);
            \draw (b1) -- (b3.apex);

            \node[vertex, label=above:{r = z}] (c1) at (9,0) {};
            \node[isosceles triangle,minimum width=1cm,minimum height=1cm,shape border rotate=180,draw] (c2) at (12,1) {$B_1$};
            \node[isosceles triangle,minimum width=1cm,minimum height=1cm,shape border rotate=180,draw] (c3) at (12,-1) {$B_2$};

            \draw (c1) -- (c2.apex);
            \draw (c1) -- (c3.apex);
        \end{tikzpicture}
    \caption{The possible branch structures of an anchor-free tree. The subtrees $B_1$ and $B_2$ must have equal height leading to the unique center of the branch. This center can be used to identify the root of $B$.}
    \label{fig:possible_branch_structures}
\end{figure}

The vertices in the center of $T$ minimize the maximum distance to any other vertex in $T$. Thus any vertex not in $Z(T)$ must have at least one vertex at distance greater than $\rad(T)$ away from it. For almost all non-path trees, we will use this fact to show that $T - Z(T)$ is an anchor.

We now restrict the possible structure of $Z(T)$. It is well known that in a tree $T$, the center $Z(T)$ is either a single vertex or a pair of adjacent vertices \cite{jordan1869assemblages} (The proof inductively removes all leaves from $T$ and proves that the center remains the same under this operation.) The following lemma shows that we cannot have anchor-free trees with a unique center vertex.

\begin{lemma}\label{lem:af_bicentered_tree}
    An anchor-free tree is bicentered.
\end{lemma}
\begin{proof}
    Suppose for contradiction that $T$ has a unique center vertex $z$. Any induced subgraph isomorphic to $T - z$ must contain $n-1$ vertices and hence must be of the form $T - v$ for some vertex $v \in V(T)$.

    By lemma~\ref{lem:af_leaf_height}, all leaves of $T$ lie at the same distance from $z$. By lemma~\ref{lem:af_root_identification}, the structure of each branch $B_z(u)$ determines its root and height uniquely. Thus the multiset of branch heights at $z$ uniquely determines the structure of $T - z$.

    Since $z$ is the unique vertex of minimum eccentricity, for every $v \neq z$ we have $\epsilon(v) > \epsilon(z)$. Equivalently, there exists a branch at $v$ whose height exceeds $\rad(T)$. In particular, the multiset of branch heights at $v$ differs from that at $z$.

    Therefore $T - v$ is not isomorphic to $T - z$ for any $v \neq z$. It follows that $T - z$ occurs uniquely as an induced subgraph of $T$, contradicting anchor-freeness. 
\end{proof}

At this point, we nearly have all the tools required to complete the proof of theorem~\ref{thm:af_trees}. The last piece we need is a symmetry result allowing us to relate the branches of removal-similar vertices in $T$.

\begin{lemma}\label{lem:af_symmetric_branches}
    Let $T$ be an anchor-free tree. For any vertex $u \in T$, the graph $T - u$ is not unique. Let $v \in T$ be such that $T - u \cong T - v$. Let $B_u(v)$ be the connected component of $T - u$ containing $v$. Then $B_u(v) \cong B_v(u)$.
\end{lemma}
\begin{proof}
    Let $\{u_1, \ldots, u_k\}$ and $\{v_1, \ldots, v_k\}$ be the neighbors of $u$ and $v$ respectively. Suppose wlog, $v \in B_u(u_1)$ and $u \in B_v(v_1)$.
    Since $T - u \cong T - v$, an isomorphism $\alpha: T - u \rightarrow T - v$ restricts to isomorphisms of the connected components of $T - u$ and $T - v$ respectively. It suffices to show that for each $j \neq 1$, $B_u(v) \ncong B_v(v_j)$. This follows since $B_u(v)$ contains $v$ and all branches of $v$ which are not $B_v(u)$ as proper subtrees.

    Thus the only branch of $T - v$ which is sufficiently large to be isomorphic to $B_u(v)$ is $B_v(u)$.
\end{proof}

\begin{figure}
    \centering
    \begin{tikzpicture}
        \node[vertex] (a1) at (0,0) {};
        \node[vertex] (a2) at (0,-1) {};
        \node[vertex] (a3) at (0,-2) {};
        \node[vertex] (b1) at (3,0) {};
        \node[vertex, label=right:z] (b2) at (3,-1) {};
        \node[vertex] (b3) at (3,-2) {};
        \node[vertex] (c1) at (6,0) {};
        \node[vertex] (c2) at (6,-1) {};
        \node[vertex] (c3) at (6,-2) {};
        \node[vertex, label=right:u] (u) at (1.5,-3) {};
        \node[vertex, label=right:v] (v) at (4.5,-3) {v};
        
        \draw (a1) -- (a2);
        \draw (a2) -- (a3);
        \draw (b1) -- (b2);
        \draw (b2) -- (b3);
        \draw (c1) -- (c2);
        \draw (c2) -- (c3);
        \draw (a2) -- (u);
        \draw (u) -- (b1);
        \draw (b2) -- (v);
        \draw (v) -- (c1);
    \end{tikzpicture}
    \caption{A tree with two pseudosimilar vertices: $T - u \cong T - v$ but $u \nsim v$. It's still true that $B_u(v) \cong B_v(u)$, but this tree still has an anchor $T - Z(T) = T - z$.}
    \label{fig:pseudosimilar_tree}
\end{figure}
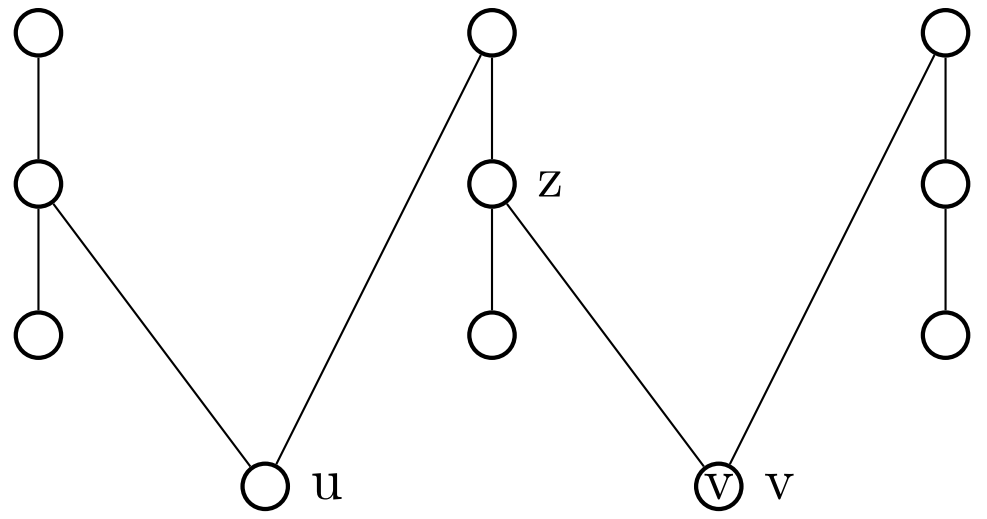

We note that lemma~\ref{lem:af_symmetric_branches} does not strictly require $T$ to adhere to a reflective, symmetric structure. For example, in Figure~\ref{fig:pseudosimilar_tree}, the branches $B_u(v)$ and $B_v(u)$ are isomorphic by a translation rather than a reflection as a rooted isomorphism. In spite of this, we note that it is precisely this lack of symmetry that ends up causing the graph of Figure~\ref{fig:pseudosimilar_tree} to not be anchor-free.

We are now ready to prove the full theorem.

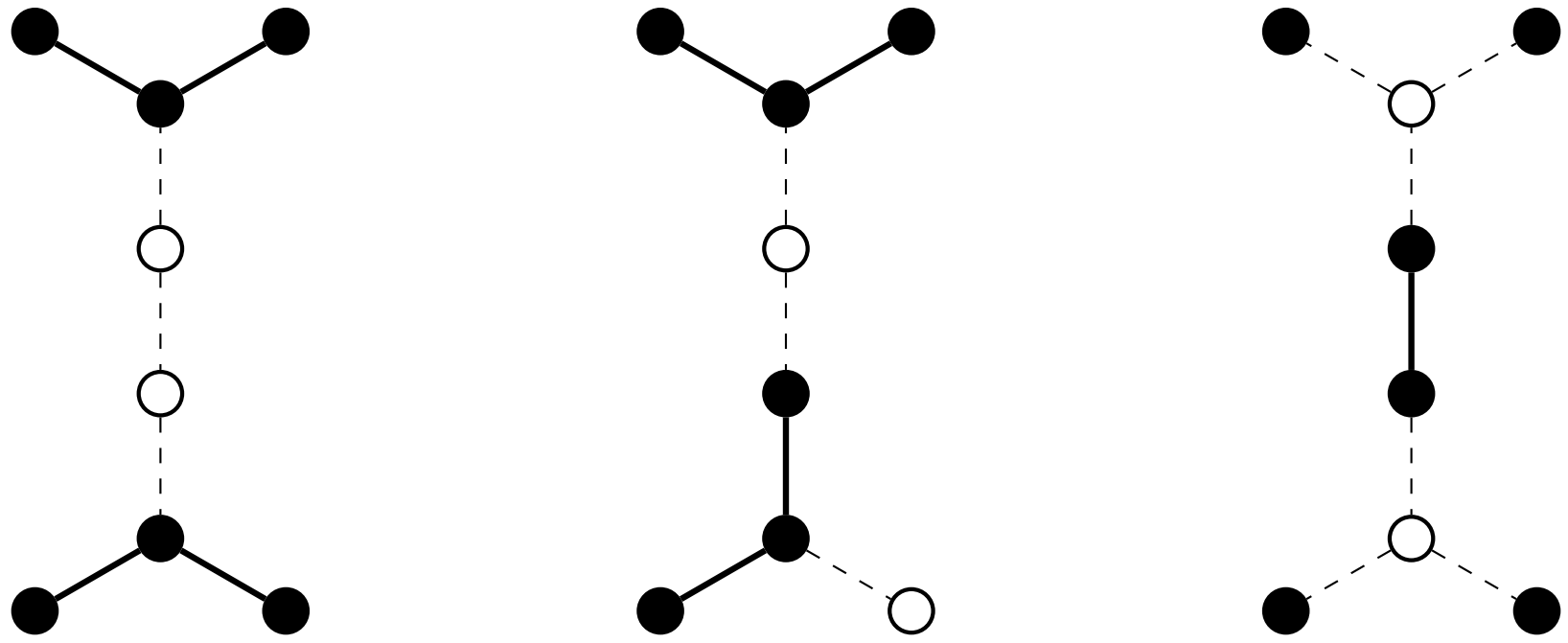
\begin{figure}
    \centering
    \def\root3{1.732}
    \begin{tikzpicture}
        \node[vertex] (z1) at (0,0) {};
        \node[vertex] (z2) at (0,-1) {};
        \node[vertex,fill=black] (r1) at (0,1) {};
        \node[vertex,fill=black] (r2) at (0,-2) {};
        \node[vertex,fill=black] (l11) at (-\root3/2,3/2) {};
        \node[vertex,fill=black] (l12) at (\root3/2,3/2) {};
        \node[vertex,fill=black] (l21) at (-\root3/2,-5/2) {};
        \node[vertex,fill=black] (l22) at (\root3/2,-5/2) {};
        
        \draw[dashed] (z1) -- (z2);
        \draw[dashed] (z1) -- (r1);
        \draw[dashed] (z2) -- (r2);
        \draw[very thick] (r1) -- (l11);
        \draw[very thick] (r1) -- (l12);
        \draw[very thick] (r2) -- (l21);
        \draw[very thick] (r2) -- (l22);
    \end{tikzpicture}
    \hspace{2cm}
    \begin{tikzpicture}
        \node[vertex] (z1) at (0,0) {};
        \node[vertex,fill=black] (z2) at (0,-1) {};
        \node[vertex,fill=black] (r1) at (0,1) {};
        \node[vertex,fill=black] (r2) at (0,-2) {};
        \node[vertex,fill=black] (l11) at (-\root3/2,3/2) {};
        \node[vertex,fill=black] (l12) at (\root3/2,3/2) {};
        \node[vertex,fill=black] (l21) at (-\root3/2,-5/2) {};
        \node[vertex] (l22) at (\root3/2,-5/2) {};
        
        \draw[dashed] (z1) -- (z2);
        \draw[dashed] (z1) -- (r1);
        \draw[very thick] (z2) -- (r2);
        \draw[very thick] (r1) -- (l11);
        \draw[very thick] (r1) -- (l12);
        \draw[very thick] (r2) -- (l21);
        \draw[dashed] (r2) -- (l22);
    \end{tikzpicture}
    \hspace{2cm}
    \begin{tikzpicture}
        \node[vertex,fill=black] (z1) at (0,0) {};
        \node[vertex,fill=black] (z2) at (0,-1) {};
        \node[vertex] (r1) at (0,1) {};
        \node[vertex] (r2) at (0,-2) {};
        \node[vertex,fill=black] (l11) at (-\root3/2,3/2) {};
        \node[vertex,fill=black] (l12) at (\root3/2,3/2) {};
        \node[vertex,fill=black] (l21) at (-\root3/2,-5/2) {};
        \node[vertex,fill=black] (l22) at (\root3/2,-5/2) {};
        
        \draw[very thick] (z1) -- (z2);
        \draw[dashed] (z1) -- (r1);
        \draw[dashed] (z2) -- (r2);
        \draw[dashed] (r1) -- (l11);
        \draw[dashed] (r1) -- (l12);
        \draw[dashed] (r2) -- (l21);
        \draw[dashed] (r2) -- (l22);
    \end{tikzpicture}
    \caption{The only non-path trees where removing the center does not yield an anchor is when each branch is a star. In this case, removing the stars' centers yields the unique anchor consisting of a single edge and many isolated vertices.}
    \label{fig:af_extended_star}
\end{figure}

\begin{theorem}
    The only anchor-free trees are the paths $P_{p-1}$ for $p$ prime.
\end{theorem}
\begin{proof}
    We have already classified the anchor-free paths as $P_{p-1}$ for $p$ prime.
    
    Let $T$ be a non-path, anchor-free tree. By lemma~\ref{lem:af_bicentered_tree}, there exists two center vertices $z_1$ and $z_2$. Since $T - z_1$ is not an anchor for $T$, there is another vertex $x$ such that $T - z_1 \cong T - x$. For any vertex $x \notin Z(T)$, the multiset of branch heights at $x$ differs from that at $z_1$ since $Z(T)$ uniquely minimizes eccentricity. Hence $T - x \cong T - z_1$ only if $x = z_2$. By lemma~\ref{lem:af_symmetric_branches}, $B_{z_1}(z_2) \cong B_{z_2}(z_1)$. Since these two disjoint branches span $T$ and are isomorphic, they must each contain exactly half of the vertices of $T$.

    Consider the subgraph $T - \{z_1, z_2\}$. Let $B$ be any branch of $T - \{z_1, z_2\}$. We have multiple cases depending on the structure of $B$.
    
    If $B$ is a single vertex, then by lemma~\ref{lem:af_leaf_height}, $\diam(T) = 3$. Since $T$ is not a path, each center vertex is adjacent to at least $2$ leaves and $T - \{z_1, z_2\}$ is a unique subgraph on $n - 2$ vertices where every connected component has exactly one vertex.

    If $B$ is a path of length at least $2$ then by lemma~\ref{lem:af_star_covers}, $T$ cannot have any vertex of degree at least $3$ and thus $T$ is a path.

    If $B$ is a star $K_{1, s}$, then either every leaf of $B$ is a leaf of $T$ or one leaf of $B$ is adjacent to the center of $T$. In the first case, every end-cutvertex of $T$ has degree $s+1$ and all leaves are at a distance $2$ from the center of $T$. This implies that all connected components of $T - \{z_1, z_2\}$ are isomorphic stars $K_{1, s}$. In this graph, removing all end-cutvertices yields an anchor consisting of $s$ isolated vertices for each end-cutvertex and the single edge $(z_1, z_2)$. (See Figure~\ref{fig:af_extended_star}.)

    In all remaining cases, the height of $B$ is at least $3$ and (up to isomorphism in the remaining star case), we can identify the root vertex of $B$ by lemma~\ref{lem:af_root_identification}. In particular, this also means that we can identify the height of $B$.

    Since $T$ is anchor-free, there exists another pair of vertices $u, v$ such that $T - \{u, v\} \cong T - \{z_1, z_2\}$. The isomorphism $T - \{z_1, z_2\} \cong T - \{u, v\}$ implies that both subgraphs have the same number of connected components, each of which has height $(D-1)/2$.
    
    Suppose without loss of generality that $u \neq z_1$ but that $u \in B_{z_2}(z_1)$ (i.e. $u$ is on the $z_1$ half of $T$). If $v$ lies on a noncentral branch of $T - u$ then the central branch $B_{u}(z_1)$ is a subgraph of $T - \{u, v\}$ with more than half the vertices of $T$ and cannot be isomorphic to any branch of $T - \{z_1, z_2\}$. Thus it must be that $v$ is on the central branch $B_u(z_1)$. Now if $v$ lies on a $u - z_1$ path, then the branch $B_v(z_1)$ is a branch of $T - \{u, v\}$ which contains greater than half of the vertices of $T$. Again, this branch cannot be isomorphic to any branch of $T - \{z_1, z_2\}$. Thus it must be that $u$ and $v$ must lie on opposite sides of the center of $T$.

    For all remaining branch structures, we can identify the height of the branches of $T - \{z_1, z_2\}$ and if $T - \{u, v\} \cong T - \{z_1, z_2\}$ then the branches must have equal heights. Since $u \neq z_1$, the only choice of $v$ that makes all branches of $T - \{u, v\}$ contained in $B_{z_1}(z_2)$ have height $(D-1)/2$ is if $v = z_2$.
    
    In the $B_{z_2}(z_1)$ half of $T$, the removal of $u$ makes all branches have height $(D-1)/2$, and thus either removing $u$ removes $z_1$ or it removes all leaves of $B_{z_2}(z_1)$. This is impossible to do by removing a single vertex of $B_{z_2}(z_1)$ since there must be at least two leaves by lemma~\ref{lem:af_star_covers}.

    Thus we conclude that for all remaining branch structures, the only choice of $\{u, v\}$ so that $T - \{u, v\} \cong T - \{z_1, z_2\}$ is $\{u, v\} = \{z_1, z_2\}$ and $T - \{z_1, z_2\}$ is an anchor for $T$.
\end{proof}

\subsection{Unique-Subgraph-Free Trees}

The previous section characterizes anchor-free trees as paths with length $p-1$ for $p$ a prime. While many of the same ideas transfer to the unique-subgraph-free case, our characterization will include a few non-path unique-subgraph-free trees.

\begin{theorem}\label{thm:usf_trees}
    Every unique-subgraph-free tree is either:
    \begin{itemize}
        \item A single vertex $K_1$,
        \item A path with prime number of vertices,
        \item A star graph, or
        \item A extended star $K_{1,s}$ where each leaf has been replaced with a copy of $K_{1,t}$ for $t \geq s - 1$.
    \end{itemize}
\end{theorem}
\addtocounter{theorem}{-1}

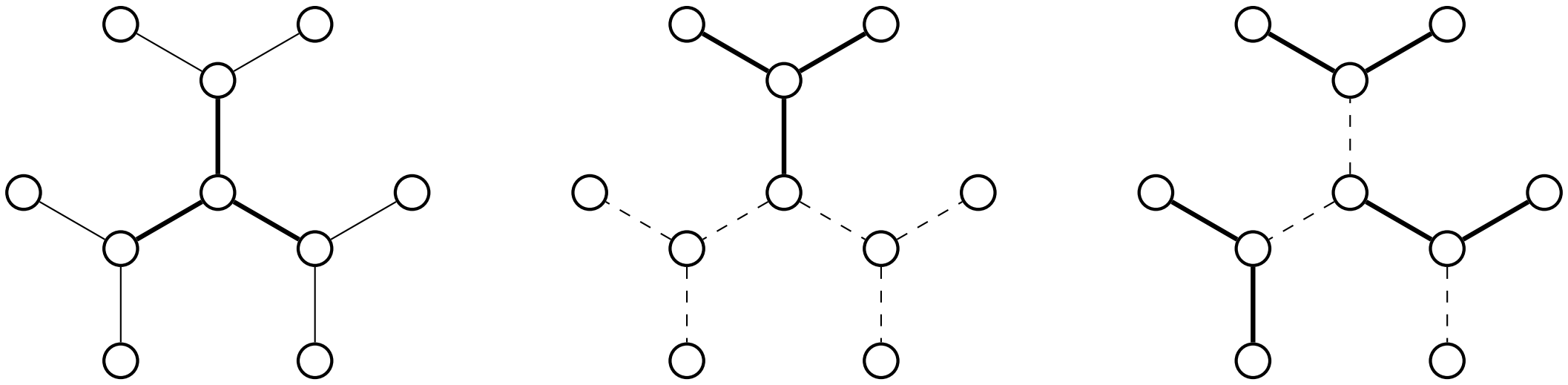
\begin{figure}
    \centering
    \def\root3{1.732}
    \begin{tikzpicture}
        \node[vertex] (z) at (0,0) {};
        \node[vertex] (r1) at (0,1) {};
        \node[vertex] (r2) at (\root3/2,-1/2) {};
        \node[vertex] (r3) at (-\root3/2,-1/2) {};

        \node[vertex] (l11) at (\root3/2,3/2) {};
        \node[vertex] (l12) at (-\root3/2,3/2) {};

        \node[vertex] (l21) at (\root3,0) {};
        \node[vertex] (l22) at (\root3/2,-3/2) {};

        \node[vertex] (l31) at (-\root3,0) {};
        \node[vertex] (l32) at (-\root3/2,-3/2) {};

        \draw[very thick] (z) -- (r1);
        \draw[very thick] (z) -- (r2);
        \draw[very thick] (z) -- (r3);

        \draw (r1) -- (l11);
        \draw (r1) -- (l12);

        \draw (r2) -- (l21);
        \draw (r2) -- (l22);

        \draw (r3) -- (l31);
        \draw (r3) -- (l32);
    \end{tikzpicture}
    \hspace{1cm}
    \begin{tikzpicture}
        \node[vertex] (z) at (0,0) {};
        \node[vertex] (r1) at (0,1) {};
        \node[vertex] (r2) at (\root3/2,-1/2) {};
        \node[vertex] (r3) at (-\root3/2,-1/2) {};

        \node[vertex] (l11) at (\root3/2,3/2) {};
        \node[vertex] (l12) at (-\root3/2,3/2) {};

        \node[vertex] (l21) at (\root3,0) {};
        \node[vertex] (l22) at (\root3/2,-3/2) {};

        \node[vertex] (l31) at (-\root3,0) {};
        \node[vertex] (l32) at (-\root3/2,-3/2) {};

        \draw[very thick] (z) -- (r1);
        \draw[dashed] (z) -- (r2);
        \draw[dashed] (z) -- (r3);

        \draw[very thick] (r1) -- (l11);
        \draw[very thick] (r1) -- (l12);

        \draw[dashed] (r2) -- (l21);
        \draw[dashed] (r2) -- (l22);

        \draw[dashed] (r3) -- (l31);
        \draw[dashed] (r3) -- (l32);
    \end{tikzpicture}
    \hspace{1cm}
    \begin{tikzpicture}
        \node[vertex] (z) at (0,0) {};
        \node[vertex] (r1) at (0,1) {};
        \node[vertex] (r2) at (\root3/2,-1/2) {};
        \node[vertex] (r3) at (-\root3/2,-1/2) {};

        \node[vertex] (l11) at (\root3/2,3/2) {};
        \node[vertex] (l12) at (-\root3/2,3/2) {};

        \node[vertex] (l21) at (\root3,0) {};
        \node[vertex] (l22) at (\root3/2,-3/2) {};

        \node[vertex] (l31) at (-\root3,0) {};
        \node[vertex] (l32) at (-\root3/2,-3/2) {};

        \draw[dashed] (z) -- (r1);
        \draw[very thick] (z) -- (r2);
        \draw[dashed] (z) -- (r3);

        \draw[very thick] (r1) -- (l11);
        \draw[very thick] (r1) -- (l12);

        \draw[very thick] (r2) -- (l21);
        \draw[dashed] (r2) -- (l22);

        \draw[very thick] (r3) -- (l31);
        \draw[very thick] (r3) -- (l32);
    \end{tikzpicture}
    \caption{The two candidate subgraphs of the extended stars are not unique.}
    \label{fig:extended_stars}
\end{figure}

It is straightforward to prove that prime paths are unique-subgraph-free using the same argument from anchor-free paths now applied to the edges of $T$. Similarly, star graphs are edge transitive and thus unique-subgraph-free by corollary~\ref{cor:transitive_properties}. The extended stars have exactly two edge orbits and thus we only need to check two different subgraphs as possible unique-subgraphs by lemma~\ref{lem:usf-respects-orbits}. See figure~\ref{fig:extended_stars} for a visual proof that these graphs are unique-subgraph-free.

As in the case of anchor-free trees, the bulk of the proof of theorem~\ref{thm:usf_trees} is left in eliminating all other possibilities. We will proceed in a similar fashion to the anchor-free case by studying the center of unique-subgraph-free trees.

\begin{lemma}\label{lem:usf_unicenter}
    With the exception of $P_2$, a unique-subgraph-free tree is unicentered.
\end{lemma}
\begin{proof}
    Suppose that a unique-subgraph-free tree is bicentered and thus has a central edge $e$. Consider the subgraph $T - e$. Any subgraph isomorphic to $T - e$ must be missing exactly one edge, i.e. $T - f \cong T - e$. Since $e$ uniquely minimizes the maximum height of the connected components of $T - e$, any other choice of edge will yield a subgraph with at least one connected component whose height is larger than both connected components of $T - e$. Thus $T - e$ is a unique subgraph.
\end{proof}

Since we now know that unique-subgraph-free trees are unicentered, let us label the center vertex at $z$ with neighbors $r_1, r_2, r_3, \ldots$. We now proceed in a similar fashion as in the proof for anchor-free trees.

\begin{lemma}\label{lem:usf_leaf_height}
    In any unique-subgraph-free tree, all leaves are at the same distance from the center.    
\end{lemma}
\begin{proof}
    Let $H$ be the subgraph which consists of all vertices and edges lying on an $\ell$-z path for each leaf $\ell$ at distance $\epsilon(z)$ from the center. By lemma~\ref{lem:usf_covering}, $H$ covers all edges of $T$. 
    
    Since there are at least $2$ such vertices at distance $\epsilon(z)$ from $z$ (otherwise $z$ would not be the center vertex), any subgraph isomorphic to $H$ in $T$ must fix $z$. If $H \neq T$, then we have a contradiction as there is no isomorphism which maps a leaf at distance $\epsilon(z)$ from $z$ to a closer leaf. Thus it must be that $H = T$ and that all leaves are at a distance $\epsilon(z)$ from $z$.
\end{proof}

As in the case of anchor-free trees, $K_{1,3}$ coverings tell us that leaves of unique-subgraph-free graphs come in at least pairs.

\begin{lemma}\label{lem:usf_star_cover}
    In any unique-subgraph-free tree $T$ with at least one vertex of degree $3$, every end-cutvertex has degree at least $3$.
\end{lemma}
\begin{proof}
    The argument is the same as in lemma~\ref{lem:af_star_covers}: $K_{1,3}$ is a subgraph and thus must cover $T$ by lemma~\ref{lem:usf_covering}. In order for $K_{1,3}$ to cover an edge incident to a leaf, the end-cutvertex must have degree at least $3$.
\end{proof}

Since the extended stars are unique-subgraph-free, the proof strategy must now deviate from anchor-free trees slightly. The next lemmas build towards proving that the extended stars are the only exceptional family added when we switch from anchor-freeness to unique-subgraph-freeness.

\begin{lemma}\label{lem:usf_center_degree}
    Let $T$ be a unique-subgraph-free graph with center $z$ and $\rad(T) \geq 2$. Then each neighbor of $z$ has equal degree. Furthermore if $\rad(T) \geq 3$, then this degree is equal to the degree of $z$.
\end{lemma}
\begin{proof}
    Consider the subgraph $T'$ obtained by removing all leaves from $T$. By lemma~\ref{lem:usf_covering}, $T'$ edge covers $T$. In particular, every leaf edge of $T$ is contained in a subgraph isomorphic to $T'$. Since $T'$ is connected and all leaves of $T$ are at the same height by lemma~\ref{lem:usf_leaf_height}, any embedding of $T'$ into $T$ which contains a leaf $\ell$ of $T$ must map $z$ one vertex closer to $\ell$. In particular, this means that the neighbor of $z$ on the $z$-$\ell$ path must have degree at least equal the degree of $z$. Since this holds for each leaf of $T$, all neighbors of $z$ must have degree at least $d(z)$.

    Now if any neighbor of $z$ has degree $c > d(z)$, then $K_{1,c}$ is a subgraph of $T$. Thus by the covering lemma again, each edge of $T$ is contained in a subgraph isomorphic to $K_{1,c}$. However, since $c > d(z)$, $z$ cannot be the center of any such subgraph. Consequently, each edge incident to $z$ must also be incident to a vertex of degree exactly $c$. Thus all neighbors of $z$ must have equal degree.

    Lastly, if $\rad(T) \geq 3$, all neighbors of $z$ have the same degree in $T'$ as they do in $T$. Now since every embedding of $T'$ into $T$ containing a leaf $\ell$ must map a neighbor of $z$ into $z$, there must be at least one neighbor of $z$ with degree equal to $z$. Since all neighbors of $z$ have equal degree, if $\rad(T) \geq 3$ then all neighbors of $z$ must have the same degree as $z$.
\end{proof}

If $\rad(T) = 1$, then $T$ is a star. Stars are edge transitive and thus unique-subgraph-free. Now if $T$ is not a path and $\rad(z) \geq 3$, then since $K_{1,3}$ is a subgraph of $T$, each edge incident to $z$ is contained in a subgraph isomorphic to $K_{1,3}$. Thus either $z$ or its neighbors have degree at least $3$ and by the previous lemma it must be that both have degree at least $3$.

We can use induction to extend the argument in the previous lemma since the embeddings mapping $T'$ into a subgraph of $T$ containing a leaf $\ell$ of $T$ maps a neighbor of $z$ into a vertex two steps closer to $\ell$ than $z$. Thus we have the following stronger version of the previous lemma.

\begin{lemma}\label{lem:usf_equal_degree}
    If $T$ is a non-path, unique-subgraph-free graph with $\rad(T) \geq 3$, then every non-leaf vertex of $T$ has the same degree.
\end{lemma}
\begin{proof}
    Let $T'$ be the subgraph obtained by removing all leaves from $T$ as before. Again let $z$ be the center of $T$. Assume for induction that all vertices at most distance $k$ from $z$ have the same degree and that the leaves of $T$ are at a distance at least $k+2$. The base case of $k = 1$ follows from lemma~\ref{lem:usf_center_degree}. Let $z - v_1 -v_2, -\ldots, - v_k - v_{k+1} - \ldots v_{r-1} - \ell$ be a path from $z$ to a leaf of $T$. Every embedding from $T'$ into $T$ which contains $\ell$ maps some vertex at distance $k$ from $z$ into $v_{k+1}$. Thus $d(v_{k+1}) \geq d(v_k) = d(z)$. Since $K_{1,d(v_{k+1})}$ is a subgraph of $T$, each edge incident to $z$ is contained in a subgraph isomorphic to $K_{1,d(v_{k+1})}$. However, this is only possible if $d(v_{k+1}) = d(z)$. Since every vertex at distance $k+1$ from $z$ lies on a path from $z$ to a leaf, the induction step holds for every vertex at distance $k+1$ from $z$.
\end{proof}

We are finally ready to prove the full theorem.

\begin{theorem}
    Every unique-subgraph-free tree is either:
    \begin{itemize}
        \item A single vertex $K_1$,
        \item A path with prime number of vertices,
        \item A star graph, or
        \item An extended star $K_{1,s}$ where each leaf has been replaced with a copy of $K_{1,t}$ for $t \geq s - 1$.
    \end{itemize}
\end{theorem}
\begin{proof}
    When $\rad(T) = 2$, lemma~\ref{lem:usf_center_degree} implies that the each neighbor of $z$ has equal degree at least $d(z)$. As long as this condition holds, the resulting graph is unique-subgraph-free since the edges partition into two orbits: adjacent to $z$ and adjacent to leaves. The copy of $K_{1, d(z)}$ in the center can be embedded into any neighbor of $z$ and the graph obtained by taking all edges adjacent to leaves is isomorphic to the same graph where any single edge is rotated from a leaf to $z$ (see Figure~\ref{fig:extended_stars}).
    
    By lemma~\ref{lem:usf_equal_degree} all remaining candidate unique-subgraph-free trees can be identified by two parameters: $d(z)$ and $\rad(T)$. Once these two parameters are chosen, we can explicitly construct $T$ by starting from the center vertex $z$.

    The only remaining cases we have left to check are trees with $\rad(T) \geq 3$. If $\rad(T) \geq 3$ and $d(z) = 2$, then $T$ must be a prime path by lemma~\ref{lem:usf_equal_degree}. 
    
    Suppose $\rad(T) \geq 3$ and $d(z) \geq 3$. In any such case, we claim that $T - z$ is a unique subgraph of $T$. Any subgraph isomorphic to $T - z$ must contain $d(z)$ isomorphic subtrees of height $\rad(T) - 1$. Since the interior vertices of $T$ all have equal degree, any subtree of $T - v$ containing $z$ must either contain all edges incident to $z$ or $z$ must be a leaf. The first case is impossible since there is not enough space to embed the remaining subtrees of $T - z$ with height $\rad(T) - 1$ while staying disconnected from the subtree containing $z$ as an internal vertex. The second case is also impossible by lemma~\ref{lem:usf_star_cover} as $z$ cannot be a singleton leaf. Thus there is no other subgraph of $T$ isomorphic to $T - z$.
\end{proof}

Combining theorems~\ref{thm:af_trees}~and~\ref{thm:usf_trees}, we see that there is exactly one tree which has neither an anchor nor a unique subgraph: $P_2$.

\section{Future Directions}
The related properties of anchor-freeness and unique-subgraph-freeness impose very strong requirements on the structure of graphs. We were able to use this structure to characterize trees with both properties. A natural next step is to extend this characterization to graphs with connectivity $\kappa(G) = 1$. For these graphs, the block-cut tree of $G$ exhibits a tree structure where it may be possible to leverage our characterization.

All known anchor-free graphs have some nontrivial symmetry. The characterization of anchor-free trees suggests that this is likely not a random observation but rather a structural requirement of anchor-free graphs. In particular, it is natural to ask whether an asymmetric anchor-free graph exists. The author is currently investigating this question. In an asymmetric anchor-free graph, the absence of automorphisms must be compensated for by numerous isomorphic copies of proper induced subgraphs; consequently, removal-similar vertices necessarily occur as pseudosimilar vertices. Understanding the structure imposed by such pseudosimilarity may therefore provide a route toward either constructing asymmetric anchor-free graphs or proving that none exist.

\section*{Acknowledgments}
The author thanks his advisors Dr. Edinah Gnang and Dr. Donniell Fishkind for helpful discussions and guidance throughout this research.

\section*{AI Usage Statement}
The usage of AI in this paper was limited to verifying the results presented as well as for editing the initial draft of the paper into a more polished and presentable version.

\bibliographystyle{alpha}
\bibliography{references}

\end{document}